\documentclass[12pt]{article}
\usepackage{amsmath,amssymb,amsfonts,epsfig,graphicx,longtable,wrapfig}

\newtheorem{theorem}{Theorem}
\newtheorem{lemma}{Lemma}

\newtheorem{corol}{Corollary}

\newcommand{\rn}{{\mathcal{L}}}

\begin{document}

\begin{center}
{\Large Zero-one laws for $k$-variable first-order logic of sparse random graphs}\\
\vspace{0.3cm}
{\large A.S. Razafimahatratra\footnote{Moscow Institute of Physics and Technology}, M. Zhukovskii\footnote{Moscow Institute of Physics and Technology, laboratory of advanced combinatorics and network applications; Adyghe State University, Caucasus mathematical center}} 

\end{center}

\section{Introduction}\label{intro}

The \emph{first-order (FO) logic} $\rn$ of undirected graphs is the set of finite \emph{first-order sentences} (see, e.g., \cite{Libkin,Survey,Strange}) from the alphabet consisting of
\begin{enumerate}
	\item symbols of variables $x,y,z,x_1,\ldots$,
	\item logical connectivities: $\wedge, \ \vee,\ \neg,\ \Rightarrow, \Leftrightarrow$,
	\item the existential and universal quantifiers: $\exists, \forall$,
	\item the relational symbols $=$ (equality) and $\sim$ (adjacency).
\end{enumerate}

We denote by $\rn^k$ the fragment of $\rn$ consisting of sentences with at most $k$ variables. For example, the sentence 
$$
\phi=\forall x\forall y \quad \biggl[\biggl(\neg[x=y]\wedge\neg[x\sim y]\biggr)\Rightarrow\biggl(\exists z\quad [z\sim x\wedge z\sim y]\biggr)\biggr]
$$ 
belongs to $\mathcal{L}^3$ and {\it expresses} the property that every two distinct non-adjacent vertices of a graph have a common neighbor, i.e., the diameter is at most $2$. 
Nevertheless, the expressive power of $\rn^3$ is much stronger than it seems to be at a first glance. In particular, for every $k$, there exists a sentence $\phi_k$ in $\rn^3$ saying that the diameter of a graph is at most $k$:
$$
 \forall x\forall y\quad\biggl[x=y\vee x\sim y \vee\biggl(\exists z\,[z\sim x\wedge z\sim y]\biggr)\vee\biggl(\exists z\,[z\sim x \wedge(\exists x\,\,[z\sim x\wedge x\sim y])]\biggr)\vee\ldots\biggr]
$$
However, the connectedness can not expressed in FO~\cite{Survey}, and it is quite natural to express it by the sentence $\phi_{\omega}$ defined in the same way as $\phi_k$ but having infinitely (countably) many disjunctions. The sentence $\phi_{\omega}$ belongs to the infinite FO logic $\rn^3_{\infty,\omega}$ (\cite{Libkin}, Chapter 11.1; the logic $\rn^k_{\infty,\omega}$ consists of all FO sentences having countably many disjunctions and conjunctions but at most $k$ variables). Recall that $\rn^{\omega}_{\infty,\omega}=\bigcup_{k=1}^{\infty}\rn^k_{\infty,\omega}$.\\

In this paper, we study asymptotical behavior of probabilities of truth of FO sentences from $\rn^k$ on the binomial random graph $G(n,p)$. Let us recall that $V_n:=\{1,\ldots,n\}$ is the set of vertices of this graph, and every edge from ${V_n\choose 2}$ appears independently from the others with probability $p$.

For an arbitrary logic $\mathcal{F}$, the random graph $G(n,p)$ is said to {\it obey the zero-one law w.r.t. $\mathcal{F}$} if, for every sentence $\phi$ from $\mathcal{F}$,   $\lim_{n\rightarrow\infty} {\sf P}(G(n,p)\models \phi) \in \{0,1\}$. If, for every sentence $\phi$ from $\mathcal{F}$, the limit $\lim_{n\rightarrow\infty} {\sf P}(G(n,p)\models \phi)$ exists (but not necessarily equal to  0 or 1), then $G(n,p)$ {\it obeys the convergence law w.r.t. to $\mathcal{F}$}.

In 1969, Glebskii et al. \cite{Glebskii} (and, independently, Fagin~\cite{Fagin} in 1976) proved that $G(n,1/2)$ obeys the zero-one law w.r.t. $\rn$ (or, simply, the FO zero-one law). In~\cite{SpencerEhren}, Spencer noticed that this result follows from a simple combinatorial argument and can be trivially extended for every $p$ such that $\min\{p,1-p\}n^{\beta}\to\infty$ for every positive constant $\beta$. In 1988, Shelah and Spencer \cite{ShelahSpencer} proved that the FO zero-one law holds for $p = n^{-\alpha}$ if and only if one of the following condition holds:
\begin{itemize}
	\item $\alpha$ is positive irrational,
	\item $\alpha >2$,
	\item $1+\frac{1}{m+1}< \alpha < 1+\frac{1}{m}$ for some $m \geq 1$.
\end{itemize}
Moreover, they proved that, for rational $\alpha\in(0,1)$, even the FO convergence law fails. In 1992, Lynch~\cite{Lynch_finite} proved that, for $\alpha=1$ and $\alpha=1+\frac{1}{m}$, the random graph $G(n,n^{-\alpha})$ obeys the FO convergence law.

In 1992, Kolaitis and Vardi~\cite{KolaitisVardi} proved that $G(n,1/2)$ obeys the zero-one law w.r.t. $\rn^{\omega}_{\infty,\omega}$. As in the finite case, a combinatorial proof works for $p$ such that $\min\{p,1-p\}n^{\alpha}\to\infty$ for every positive constant $\alpha$ as well. 

Surely, from the result of Spencer and Shelah, it follows that $G(n,n^{-\alpha})$ does not obey the convergence law w.r.t. $\mathcal{L}^{\omega}_{\infty,\omega}$ for rational $\alpha\in(0,1)$. For $\alpha>2$, the random graph obeys the zero-one law w.r.t.  $\mathcal{L}^{\omega}_{\infty,\omega}$ since a.a.s. (with asymptotical probability $1$) this graph is empty. For $\alpha\in(1,2]$, Lynch in 1993~\cite{Lynch_infinite} proved that the infinite case mirrors the finite one: there is no zero-one law w.r.t. $\mathcal{L}^{\omega}_{\infty,\omega}$ if and only if $\alpha=1+\frac{1}{m}$ for some positive integer $m$. And, for $\alpha=1+\frac{1}{m}$, $G(n,n^{-\alpha})$ obeys the convergence law w.r.t. $\mathcal{L}^{\omega}_{\infty,\omega}$. In 1993, Tyszkiewicz~\cite{Tyszk_infinite} proved that $G(n,\frac{1}{n})$ does not obey the convergence law w.r.t. $\mathcal{L}^{\omega}_{\infty,\omega}$. The last case $\alpha\in(0,1)\setminus\mathbb{Q}$ was solved only recently by Shelah~\cite{Shelah_preprint}. In~\cite{McArthur}, McArthur claimed that there is no convergence law w.r.t. $\mathcal{L}^{\omega}_{\infty,\omega}$ for such $\alpha$ referring to a joint paper with Spencer that never appeared. In 2017, Shelah~\cite{Shelah_preprint} proved that the claim is true.\\

The aforementioned results imply that, for every $\alpha\in(0,1)$, there exists $k$ such that $G(n,n^{-\alpha})$ does not obey the convergence law w.r.t. $\rn^k_{\infty,\omega}$. But how large should such $k$ be? In 1997, McArthur~\cite{McArthur} proved that, for $\alpha<\frac{1}{k-1}$, $G(n,n^{-\alpha})$ obeys the zero-one law w.r.t. $\rn^k_{\infty,\omega}$, whereas $G(n,n^{-1/(k-1)})$ does not obey even the convergence law w.r.t. $\mathcal{L}^k_{\infty,\omega}$ for $k\geq 4$. For $k\in \{2,3\}$, $G(n,n^{-1/(k-1)})$ obeys the zero-one law w.r.t. $\rn^k_{\infty,\omega}$.\\

Surely, from the result of McArthur, it follows that, for $k\geq 4$, $G(n,n^{-\alpha})$ obeys the zero-one law w.r.t. $\mathcal{L}^k$ for $\alpha<\frac{1}{k-1}$. But what about $\alpha=\frac{1}{k-1}$? If the zero-one law still holds in this case, then is there any $\varepsilon>0$ such that the zero-one law holds for every $\alpha\in(\frac{1}{k-1},\frac{1}{k-1}+\varepsilon)$? In this paper, we give the following answer:

\begin{theorem}
For every $k\geq 2$, the random graph $G(n,n^{-1/(k-1)})$ obeys the zero-one law w.r.t. $\mathcal{L}^k$. Moreover, for every $k\geq 3$ and every $\varepsilon>0$, there exists $\alpha\in(\frac{1}{k-1},\frac{1}{k-1}+\varepsilon)$ such that $G(n,n^{-\alpha})$ does not obey the zero-one law w.r.t. $\mathcal{L}^k$.
\label{thm}
\end{theorem}

To the best of our knowledge, the negative result of Theorem~\ref{thm} was not known before even for $\mathcal{L}^k_{\infty,\omega}$. Moreover, we have not found any known result on validity of the zero-one law w.r.t. $\mathcal{L}^k$ or $\mathcal{L}^k_{\infty,\omega}$ for $\alpha>\frac{1}{k-1}+\varepsilon$ and some positive $\varepsilon$. For large $\alpha$, in the context of zero-one laws, something is known for the smaller class $\mathcal{L}_k$ of sentences having quantifier rank (this notion is defined in, e.g., \cite[Definition 3.8]{Libkin}) at most $k$, see~\cite{OstrZhukIzv,Zhuk_1,Zhuk_0}.

As we mentioned above, for every irrational $\alpha\in(0,1)$ and large enough $k$, the random graph $G(n,n^{-\alpha})$ does not obey the zero-one law w.r.t. $\mathcal{L}^k_{\infty,\omega}$. But are there, for large enough $k$, such irrational $\alpha$ close to $\frac{1}{k-1}$? Moreover, what is the minimum $k$ such that, for some irrational $\alpha\in(0,1)$, the random graph $G(n,n^{-\alpha})$ does not obey the zero-one law w.r.t. $\mathcal{L}^k_{\infty,\omega}$? We do not have answers on both questions.

We also tried to prove or disprove the convergence near $\frac{1}{k-1}$, but we failed. So, we leave it as an open question: is there any $\varepsilon>0$ such that the convergence law w.r.t. $\mathcal{L}^k_{\infty,\omega}$ holds for every $\alpha\in(\frac{1}{k-1},\frac{1}{k-1}+\varepsilon)$?\\ 

This work is organized in the following way.  In Section~\ref{prem}, we review the tools we use in our proof. In Sections~\ref{proof1}~and~\ref{proof2}, we prove Theorem~\ref{thm}.

\section{Preliminaries}\label{prem}	

\subsection{Extensions}
\label{ext_prem}

Everywhere below in this section, $G$ and $H$ are two graphs such that $H$ is a subgraph of $G$ ($H\subset G$);  the number of vertices of $H$ equals $v(H)=\ell$, and $v(G)-v(H)=m$. In particular, $u_1,\ldots,u_{\ell}$ are the vertices of $H$, and $v_1,\ldots,v_m$ are the vertices of $G$ outside $H$.


Let $\alpha>0$. Below, we denote by $e(G)$ the number of edges of $G$. In the paper, we use the following notions that can be also found, e.g., in~\cite{AlonSpencer}, Chapter 10:
	\begin{enumerate}
		\item The pair $(G,H)$ is $\alpha$-\emph{sparse} if $\left[ v(G)-v(H) \right] -\alpha\left[ e(G) - e(H)\right]>0$.
		\item The pair $(G,H)$ is $\alpha$-\emph{dense} if $\left[ v(G)-v(H) \right] -\alpha\left[ e(G) - e(H)\right]<0$.
		\item The pair $(G,H)$ is $\alpha$-\emph{safe} if for any subgraph $S$ such that $G\subset S \subseteq H$, $(G,S)$ is $\alpha$-sparse.
	\end{enumerate}
	
We will frequently use the following result 
(see, e.g., \cite{AlonSpencer}, Chapter 10; \cite{Spencer_extensions_counting}, \cite{Spencer_extensions_threshold}).


\begin{lemma}
Let $\alpha = \frac{1}{k-1}$ and $(G,H)$ be $\alpha$-safe. Then, with asymptotical probability $1$, for any distinct vertices $x_1,\ldots.x_{\ell}$ of $G(n,n^{-\alpha})$, there exist distinct vertices $y_1,\ldots,y_m$ such that 
\begin{equation}
 \left(\bigwedge_{i,j}\left(x_i\sim y_j\Leftrightarrow u_i\sim v_j\right)\right)\wedge
\left(\bigwedge_{i\neq j}\left(y_i\sim y_j\Leftrightarrow v_i\sim v_j \right)\right)
\label{extension_property}
\end{equation}
 and any $k-1$ vertices among $x_1,\ldots,x_{\ell},y_1,\ldots,y_m$ such that at least one of them belongs to $\{y_1,\ldots,y_m\}$ have no common neighbors outside $\{x_1,\ldots,x_{\ell},y_1,\ldots,y_m\}$.
\label{generic_extensions}
\end{lemma}



\subsection{Subgraph containment}\label{subgraphs}

Given a graph $H$, let $\rho(H)=\frac{e(H)}{v(H)}$, and $\rho^{\max}(H)=\max_{\tilde H\subseteq H}\rho(\tilde H)$. A graph $H$ is called {\it strictly balanced}, if, for every its proper subgraph $\tilde H$, $\rho(\tilde H)<\rho(H)$. We will use the following well-known results on asymptotical probabilities of the property of containing subgraphs isomorphic to $H$.

\begin{theorem}[A. Ruci\'{n}ski, A. Vince, 1986 \cite{RV}]
Let $H$ be an arbitrary non-empty graph. Then, for $p\gg n^{-1/\rho^{\max}(H)}$, ${\sf P}(G(n,p)\supset H)\to 1$ as $n\to\infty$. If $p\ll n^{-1/\rho^{\max}(H)}$, then ${\sf P}(G(n,p)\supset H)\to 0$ as $n\to\infty$.
\label{small_thr}
\end{theorem}

\begin{theorem}[B. Bollob\'{a}s, 1981 \cite{Bol_small}]
Let $H$ be a strictly balanced graph, $c>0$, $p=cn^{-1/\rho(H)}$. Then the number of subgraphs in $G(n,p)$ isomorphic to $H$ converges in distribution to a Poisson random variable with the parameter $c^{e(H)}/a(H)$, where $a(H)$ is the number of automorphisms of $H$.
\label{pois}
\end{theorem}

\subsection{The $k$-Pebble game}\label{k-pebble}

The maximum number of nested quantifiers of a first-order sentence is called its {\it quantifier rank} (see the formal definition in \cite[Definition 3.8]{Libkin}). Denote $\mathcal{L}_N$ the set of all first-order sentences having quantifier rank at most $N$.

Let $G$ and $H$ be two graphs and let $k$ be a positive integer. Let us recall the rules of the $k$-\emph{pebble Ehrenfeucht-Fra{\"\i}ss{\'e} game} on $G$ and $H$ played in $N$ rounds~\cite[Chapter 11.2]{Libkin}.

Two players, \emph{Spoiler} and \emph{Duplicator} (or \emph{he} and \emph{she} resp.) have equal sets of $k$ pairwise different pebbles. In each round, Spoiler takes a pebble and puts it on a vertex in $G$ or in $H$;  then Duplicator has to put her copy of this pebble on a vertex of the other graph. Note that the pebbles can be reused and change their positions during the play. Duplicator's objective is to ensure that the pebbling determines a partial isomorphism between $G$ and $H$ after each round up to $N$; when she fails, she immediately loses.

We will use the following Ehrenfeucht-Fraiss\'e type theorem for the pebble game (the proof can be found in~\cite{Immerman-book}).
\begin{theorem}
Let $G$ and $H$ be two graphs and let $N$ be a positive integer. Then Duplicator has a winning strategy in $N$ rounds of 	the $k$-pebble game on $G,H$ if and only if there is no sentence from $\rn^{k}\cap \rn_N$ that distinguishes between $G$ and $H$ (that is, no sentence $\phi$  such that $G \models \phi$ and $H\not\models \phi$).
\label{Ehren}
\end{theorem}

It is well-known that elementary equivalence results as above imply necessary and sufficient conditions for validity of zero-one laws. For the logic $\mathcal{L}_N$, the corollary from the relevant Ehrenfeucht-Fraiss\'e type theorem is stated and proven in, e.g., \cite[Theorem 7]{Survey}, \cite[Theorem 2.5.1]{Strange}. In our case, the arguments are the same. So, we only present the corollary from Theorem~\ref{Ehren} and omit its proof. 

\begin{corol}
	Let $k$ be a fixed positive integer. The random graph $G(n,p)$ obeys the zero-one law w.r.t. $\rn^k$ if and only if, for every positive integer $N$, 
	\begin{align*}
	\displaystyle \mathsf{P}\left[
	\begin{aligned}
	\mbox{Du}&\mbox{plicator has a winning strategy in }\\ &\mbox{the }k\mbox{-pebble game on }G(n,p),\,G(m,p)\mbox{ in }N\mbox{ rounds} 
	\end{aligned}
	\right] \xrightarrow[n,m \rightarrow\infty]{}1,
	\end{align*}\label{zero_one_law_EF}
where the random graphs $G(n,p)$ and $G(m,p)$ are independent.
\label{corol}
\end{corol}

\section{The proof of the first part of Theorem~\ref{thm}}\label{proof1}

Here, we prove that, for every $k\geq 2$, the random graph $G(n,n^{-1/(k-1)})$ obeys the zero-one law w.r.t. $\mathcal{L}^k$. For $k\in\{2,3\}$ this follows from the McArthur's result. So, we consider $k\geq 4$.

Our main tool is Corollary~\ref{corol}. In other words, we will prove that, for every $N$, a.a.s. (asymptotically almost surely) Duplicator wins the $k$-pebble game on $G(n,n^{-1/(k-1)})$ and $G(m,m^{-1/(k-1)})$ in $N$ rounds.\\

Let us first notice that existence of a winning strategy of Duplicator in the $(k-1)$-pebble game clearly follows from the, so called, {\it $(k-1)$-extension property} (i.e., for every sets of vertices $U\subseteq W$ such that $|W|<k-1$, there exists a common neighbor of vertices of $U$ having no neighbors among the vertices of $W\setminus U$): if both graphs have the $(k-1)$-extension property, then Duplicator has a winning strategy in the $(k-1)$-pebble game in arbitrarily large number of rounds. It is very well known (see, e.g., the proof of Theorem 22 in~\cite{Survey}) that a.a.s. $G(n,n^{-1/(k-1)})$ has the $(k-1)$-extension property which implies the zero-one law w.r.t. $\mathcal{L}^{k-1}$. However, for $k$ pebbles, this method does not work, since, by Lemma~\ref{generic_extensions}, a.a.s. there exist $(k-1)$-sets having no common neighbors.\\ 


Let $N$ be a positive integer and $\mu=\mu(N)$ be as large as desired. Let graphs $G_1$ and $G_2$ have the following properties:
\begin{itemize}
\item[1)] $(k-1)$-extension property;
\item[2)] the property from Lemma~\ref{generic_extensions}: for every $\frac{1}{k-1}$-sparse pair $(G,H)$ (notations are from Section~\ref{ext_prem}) such that $\ell\leq k-1$, $m\leq \mu$, and any distinct vertices $x_1,\ldots.x_{\ell}$ there exist distinct vertices $y_1,\ldots,y_m$ such that (\ref{extension_property}) holds,  and any $k-1$ vertices among $x_1,\ldots,x_{\ell},y_1,\ldots,y_m$ such that at least one of them belongs to $\{y_1,\ldots,y_m\}$ have no common neighbors outside $\{x_1,\ldots,x_{\ell},y_1,\ldots,y_m\}$.
\end{itemize}

From the above arguments and Lemma~\ref{generic_extensions}, it is enough to prove that Duplicator wins the $k$-pebble game on graphs $G_1$ and $G_2$ in $N$ rounds. Without loss of generality, we assume that Spoiler and Duplicator place all $k$ pairs of pebbles in the first $k$ rounds.

Let us enumerate the pairs of pebbles by $1,\ldots,k$ and denote by $a_r(i)$ and $b_r(i)$ the vertices pebbled by the $i$-th pebbles after the $r$-th round in the graphs $G_1$ and $G_2$ respectively.\\

We skip the first $k-2$ rounds, since the trivial strategy of Duplicator follows from the $(k-1)$-extension property.\\

Let us first describe the main ingredient of the winning strategy of Duplicator. Clearly, if, in every round from $k$ up to $N$, Spoiler chooses a `good' vertex which is not a common neighbor of the other pebbled vertices, than Duplicator wins by the property 2). Unfortunately, Spoiler is allowed to choose common neighbors. Nevertheless, Duplicator may take a look on all such possible future moves of Spoiler in advance and consider the graph induced by all these moves. Below, we represent all these moves in the, so called, strategy tree $T^{\infty}$. This induced graph observed by Duplicator, together with the currently pebbled subgraph of Spoiler's graph, would give an $\frac{1}{k-1}$-sparse pair, and so, at the first glance, Duplicator may find a copy of this graph that extends the pebbled subgraph of the graph she currently plays in. But the problem is that the size of this graph could be arbitrarily large (in particular, much bigger than $\mu$). Fortunately, she may remove from the observed `extension' a significant part in a way such that the obtained graph would still contain all possible different `bad' moves of Spoiler. Below, these remaining moves are represented in the, so called, refinement of $T^{\infty}$. Let us formalize this.\\

Assume that $r\leq N-1$ rounds are played, and, either $r=k-2$, or, in the round $r+1$, Spoiler moves a pebble to a vertex which is not a common neighbor of the previously pebbled vertices. In the latter case, without loss of generality, we assume that Spoiler moves the pebble $k$. In both cases, without loss of generality, we assume that, in the round $r+1$, Spoiler makes a move in $G_1$.\\
 

Given a graph $G$, consider a rooted tree $T$ such that the non-root vertices of $T$ are triplets $x=(R=R(x),S=S(x),v=v(x))$, where $R\in{\{1,\ldots,k\}\choose k-1}$, $S\in [V(G)]^{k-1}$ and $v\in V(G)$ is a common neighbor of the vertices from $S$ in $G$. We call $T$ a {\it strategy tree of $G$}. We use this notion since the strategy trees defined below represent possible strategies of Spoiler moving pebbles only in $G$. More precisely, each of them describes all possible moves to common neighbors of the other pebbled vertices: move of a pebble $\beta$ in a round $i+r+1$, $i\in\{1,\ldots,N-r-1\}$ to a common neighbor $v$ of pebbled vertices $a_{i+r}(j)$, $j\in\{1,\ldots,k\}\setminus\{\beta\}$, is represented by the vertex $(\{1,\ldots,k\}\setminus\{\beta\},(a_{i+r}(j), j\in\{1,\ldots,k\}\setminus\{\beta\}),v)$ of the strategy tree. So, the first element of a vertex of a strategy tree represents a set of non-moved pebbles, the second element represents the set of vertices pebbled by them, and the last element represents a currently pebbled vertex.  The formal constructions are given below (after few more definitions related to general strategy trees).

Two strategy trees $T_1,T_2$ (not necessarily of the same graph $G$) are {\it $\cong_1$-isomorphic}, if there exists an isomorphism $f:\,T_1\to T_2$ of rooted trees such that, for every non-root vertex $x$ of $T_1$, its first element equals to the first element of its image: $R(x)=R(f(x))$. 

Let $x$ be a vertex of a rooted tree $T$. Then $T_x$ is the subtree of $T$ rooted in $x$ and induced on the set $\{$all the successors of $x\}\cup\{x\}$ ($y$ is a {\it successor} of $x$,  if there is a path between $R$ and $y$ having $x$ as an inner vertex). We call such a tree {\it pendant subtree of $T$}. Clearly, a pendant subtree of a strategy tree is a strategy tree as well.\\

Let us construct a strategy tree $T^{\infty}=T^{\infty}[G_1;\,{\bf a}]$ (where ${\bf a}$ is either the $(k-1)$-tuple, or $k$-tuple of currently pebbled in $G_1$ vertices) rooted in a root $\mathcal{R}_0$ recursively by adding, at step $i$, all the vertices at the distance $i$ from the root. Notice that the superscript $\infty$ in $T^{\infty}$ means that its size is not bounded by any constant rather than infiniteness of the tree.

For a set of vertices $S$, we denote by $N(S)$ the set of all its common neighbors (it is always clear from the context, which graph is considered)
.\\

If $r=k-2$, then the children of $\mathcal{R}_0$ in $T^{\infty}$ are the triplets
$$
\mathcal{R}_{1,j}=\biggl(R_{1,j}=\{1,\ldots,k-1\},\,S_{1,j}=(a_{k-1}(1),\ldots,a_{k-1}(k-1)),\,v_{1,j}\biggr),
\quad v_{1,j}\in N(\{S_{1,j}\}).
$$
If $r\geq k-1$, then the children of $\mathcal{R}_0$ in $T^{\infty}$ are the triplets
$$
\mathcal{R}_{1,j}=\biggl(R_{1,j}=\{1,\ldots,k\}\setminus\{\beta\},\,S_{1,j}=(a_{r+1}(1),\ldots,a_{r+1}(\beta-1),a_{r+1}(\beta+1),\ldots,a_{r+1}(k)),\, v_{1,j}\biggr),
$$
$$
\beta\in\{1,\ldots,k\}, v_{1,j}\in N(\{S_{1,j}\}).
$$
Suppose that we have constructed all vertices $\mathcal{R}_{i,j}=(R_{i,j},S_{i,j},v_{i,j})$ of $T^{\infty}$ at the distance $i<N-r-1$ from $\mathcal{R}_0$, where $R_{i,j}$ are $(k-1)$-subsets of $\{1,\ldots,k\}$, $S_{i,j}=(s_{i,j}^1,\ldots,s_{i,j}^{k-1})$ are $(k-1)$-tuples of vertices of $G_1$, and $v_{i,j}=:s_{i,j}^k$ are vertices of $N(\{S_{i,j}\})$. Then, the children of $\mathcal{R}_{i,j}$ are all the possible triplets $(R,S,v)$ (if they exist --- otherwise, $\mathcal{R}_{i,j}$ becomes a leaf of $T^{\infty}$), where, for some $\beta\in\{1,\ldots,k\}$, $R=\{1,\ldots,k\}\setminus\{\beta\}$, $S=(s_{i,j}^{\gamma},\,\gamma=1,\ldots,\beta-1,\beta+1,\ldots,k)$, and $v\in N(\{S\})$.

The construction finishes at step $i=N-r-1$.

Let $d\leq N-r-1$ be {\it the depth} of $T^{\infty}$ (i.e., the number of edges in a longest path starting in the root).\\

Clearly, the size of $T^{\infty}$ is not bounded by any constant. However, it has many `equivalent' pendant subtrees, and, therefore, admits a bounded `refinement'. 

More formally, 
for every $i=1,2,\ldots,d-1$ (recursively in the ascending order), and every vertex $x\in V(T^{\infty})$, consider the set of all pendant subtrees of $T^{\infty}$ rooted in children of $x$ and having depth $i$. We remove from every $\cong_1$-isomorphism class on this set all but one subtrees from $T^{\infty}$. At the end, we get the tree $T$. Clearly, the size of $T$ is bounded from above by a function of $k$ and $N$. {\it Let the parameter $\mu$ be at least this bound}. We call $T$ {\it a refinement of $T^{\infty}$}.\\

The main goal of Duplicator is to pebble a vertex $b_{k-1}(k-1)$ (here, we consider the case $r=k-2$; in the second case, Duplicator should pebble $b_{r+1}(k)$) such that 
$$
T^{\infty}[G_2;\,(b_{k-1}(1),\ldots,b_{k-1}(k-1))]\cong_1 T\quad (T^{\infty}[G_2;\,(b_{r+1}(1),\ldots,b_{r+1}(k))]\cong_1 T).
$$
Let us show that this is possible.


Let us construct a graph $A$ which is defined by the tree $T$ in the following way. Start from the graph $A_1=G_1|_{\{S_{1,j}\}}$ (notice that sets $S_{1,j}$, in fact, do not depend on $j$). For $i\in\{1,\ldots,d\}$, assume that the graph $A_i$ is constructed. The graph $A_{i+1}$ is obtained from $A_i$ by adding, for every $\mathcal{R}_{i,j}$, the vertex $v_{i,j}$ adjacent to all the vertices of $S_{i,j}$ (that belong to $A_i$ by induction) and non-adjacent to all the other vertices of $A_i$. Set $A=A_{d+1}$.

Let $B=A|_{\{a_{k-1}(1),\ldots,a_{k-1}(k-2)\}}$ if $r=k-2$, and $B=A|_{\{a_{r+1}(1),\ldots,a_{r+1}(k-1)\}}$, otherwise.

Clearly, the pair $(A,B)$ is $1/(k-1)$-safe. By the second property of $G_2$, if $r=k-2$, there exist an induced subgraph $\tilde A$ of $G_2$ such that 

\begin{itemize}

\item $b_{k-2}(1),\ldots,b_{k-2}(k-2)\in V(\tilde A)$,

\item there exists an isomorphism $\Phi:\,A\to\tilde A$ sending $a_{k-1}(i)$ to $b_{k-2}(i)$, $i\in\{1,\ldots,k-2\}$,
 
\item any $k-1$ vertices of $\tilde A$ do not have a common neighbor outside $\tilde A$.
 
\end{itemize}
 
Duplicator pebbles $b_{k-1}(k-1)=\Phi(a_{k-1}(k-1))$. Clearly, $T^{\infty}[G_2;\,b_{k-1}(1),\ldots,b_{k-1}(k-1)]\cong_1 T$ (for the sake of convenience, here and below we assume that the $\cong_1$-isomorphism sends $\mathcal{R}_{i,j}$ of $T^{\infty}[G_2;\,b_{k-1}(1),\ldots,b_{k-1}(k-1)]$ to $\mathcal{R}_{i,j}$ of $T$ for all $i$ and $j$).

If $r\geq k-1$, then there exist an induced subgraph $\tilde A$ of $G_2$ such that

\begin{itemize}

\item $b_r(1),\ldots,b_r(k-1)\in V(\tilde A)$,

\item there exists an isomorphism $\Phi:\,A\to\tilde A$ sending $a_{r+1}(i)$ to $b_r(i)$, $i\in\{1,\ldots,k-1\}$,

\item any $k-1$ vertices of $\tilde A$ having at least one vertex outside $\{b_r(1),\ldots,b_r(k-1)\}$ do not have a common neighbor outside $\tilde A$.

\end{itemize}

Duplicator moves the pebble $k$ to $\Phi(a_{r+1}(k))$. Clearly, $T^{\infty}[G_2;\,b_{r+1}(1),\ldots,b_{r+1}(k)]\cong_1 T$.\\

Finally, assume that $k-1\leq r\leq N-1$ rounds are played; in the round $r+1$, Spoiler moves a pebble to a vertex which is a common neighbor of the previously pebbled vertices. Let $\tilde r\geq k-2$ be the minimum number such that in all the rounds $\tilde r+1,\tilde r+2,\ldots,r+1$ Spoiler pebbled common neighbors of $k-1$ previously pebbled vertices. In the $\tilde r$-th round, the strategy trees $T^{\infty}[G_1;\,\mathbf{a}]$ rooted in $\mathcal{R}_0^1$ and $T^{\infty}[G_2;\,\mathbf{b}]$ rooted in $\mathcal{R}_0^2$ were constructed, where $\mathbf{a}$ and $\mathbf{b}$ are tuples of vertices pebbled in $G_1$ and $G_2$ respectively after the round $\tilde r$. Without loss of generality, assume that in the round $\tilde r$, Spoiler's move was in $G_1$. In rounds $\tilde r+1,\tilde r+2,\ldots,r$, in $G_1$ and $G_2$ the pebbles $g_{\tilde r+1},\ldots,g_r$ were moved respectively, and there exists a refinement $T$ of $T^{\infty}[G_1;\,\mathbf{a}]$ such that 
\begin{itemize}
\item $T^{\infty}[G_2;\,\mathbf{b}]\cong_1 T$,
\item in $T^{\infty}[G_2;\,\mathbf{b}]$ and $T$, there are paths $\mathcal{R}_0^2\mathcal{R}_{1,1}^2\ldots\mathcal{R}_{r-\tilde r,1}^2$ and $\mathcal{R}_0^1\mathcal{R}_{1,1}^1\ldots\mathcal{R}_{r-\tilde r,1}^1$ respectively having 
$$
v(\mathcal{R}_{i,1}^1)=a_{\tilde r+i}(g_{\tilde r+i}),\quad v(\mathcal{R}_{i,1}^2)=b_{\tilde r+i}(g_{\tilde r+i}),
$$
$$
R(\mathcal{R}_{i,1}^1)=R(\mathcal{R}_{i,1}^2)=(1,\ldots,g_{\tilde r+i}-1,g_{\tilde r+i}+1,\ldots,k),
$$ 
$i\in\{1,\ldots,r-\tilde r\}$.
\end{itemize}

Without loss of generality, assume that, in the current $(r+1)$-th round, Spoiler moves the $k$-th pebble. 

If he makes it in $G_1$, then, in $T^{\infty}[G_1;\,\mathbf{a}]$, there is a child $\mathcal{R}_{r-\tilde r+1,j}^1$ of $\mathcal{R}_{r-\tilde r,1}^1\in V(T)$ such that $v(\mathcal{R}_{r-\tilde r+1,j}^1)=a_{r+1}(k)$, and $R(\mathcal{R}_{r-\tilde r+1,j}^1)=(1,\ldots,k-1)$. Since $T$ is a refinement of $T^{\infty}[G_1;\,\mathbf{a}]$, there is child, say, $\mathcal{R}_{r-\tilde r+1,1}^1$ of $\mathcal{R}_{r-\tilde r,1}^1$ in $T$ such that $T_{\mathcal{R}_{r-\tilde r+1,1}^1}$ is $\cong_1$-isomorphic to a refinement of $T^{\infty}[G_1;\,\mathbf{a}]_{\mathcal{R}_{r-\tilde r+1,j}^1}$. Consider the (only) $\cong_1$-isomorphism $\Phi:\,T\to T^{\infty}[G_2;\,\mathbf{b}]$ and let $\mathcal{R}_{r-\tilde r+1,1}^2=\Phi(\mathcal{R}_{r-\tilde r+1,1}^1)$. Duplicator moves the pebble to $v(\mathcal{R}_{r-\tilde r+1,1}^2)$. To finish the inductive argument, it is sufficient to replace $T$ with a refinement of $T^{\infty}[G_1;\,\mathbf{a}]$ that contains the path $\mathcal{R}_0^1\mathcal{R}_{1,1}^1\ldots\mathcal{R}_{r-\tilde r,1}^1\mathcal{R}_{r-\tilde r+1,j}^1$.

If Spoiler moves in $G_2$, then, in $T^{\infty}[G_2;\,\mathbf{b}]$, there is a child, say, $\mathcal{R}_{r-\tilde r+1,1}^2$ of $\mathcal{R}_{r-\tilde r,1}^2$ such that $v(\mathcal{R}_{r-\tilde r+1,1}^2)=b_{r+1}(k)$, and $R(\mathcal{R}_{r-\tilde r+1,1}^2)=(1,\ldots,k-1)$. Consider the $\cong_1$-isomorphism $\Phi:\,T^{\infty}[G_2;\,\mathbf{b}]\to T$ and let $\mathcal{R}_{r-\tilde r+1,1}^1=\Phi(\mathcal{R}_{r-\tilde r+1,1}^2)$. Duplicator moves the pebble $k$ to $v(\mathcal{R}_{r-\tilde r+1,1}^1)$.\\

The strategy is winning for Duplicator, this finishes the proof.

\section{The proof of the second part of Theorem~\ref{thm}}\label{proof2}

Here, for every $\varepsilon>0$, we find an $\alpha\in(\frac{1}{k-1},\frac{1}{k-1}+\varepsilon)$ such that $G(n,n^{-\alpha})$ does not obey the zero-one law w.r.t. $\mathcal{L}^k$. For this, we use the following construction from~\cite{McArthur}. Consider $k-2$ non-adjacent vertices $v_1,\ldots,v_{k-2}$, and vertices $x_0,x_1,\ldots,x_{\ell},x$ with the following edges: 
\begin{itemize}
\item $x_0$ is adjacent to $v_1,\ldots,v_{k-2}$;
\item for every $i\in\{1,\ldots,\ell\}$, $x_i$ is adjacent to $x_{i-1}$ and each of $v_1,\ldots,v_{k-2}$;
\item $x$ is adjacent to $v_1,\ldots,v_{k-3},x_{\ell-1},x_{\ell}$.
\end{itemize}
The obtained graph is called {\it a terminated $(k-1)$-chain rooted in $v_1,\ldots,v_{k-2}$ of length $\ell$}. We denote it by $H^{k-1}_{\ell}$. Let us prove the following property of terminated $(k-1)$-chains.

\begin{lemma}
For every $k\geq 4$, there exists $\ell_0$ such that, for every integer $\ell\geq\ell_0$, the graph $H^{k-1}_{\ell}$ is strictly balanced.
\label{H_str_b}
\end{lemma}

{\it Proof} The graph $H^{k-1}_{\ell}$ can be constructed in the following way: start with $k-2$ non-adjacent vertices $v_1,\ldots,v_{k-2}$ and their common neighbor;  at every step $1\leq i\leq \ell$, add a vertex adjacent to all $v_1,\ldots,v_{k-2}$ and one of the previously introduced vertices; at the last step $i=\ell+1$, add a vertex adjacent to $v_1,\ldots,v_{k-3}$ and two of the previously introduced vertices. Therefore, for a proper $\tilde H\subset H^{k-1}_{\ell}$, there exist $a\leq k-2$ (the number of vertices of $\tilde H$ among $v_1,\ldots,v_{k-2}$) and $b\leq\ell+2$ (the number of vertices of $H$ among $x_0,x_1,\ldots,x_{\ell},x$) (one of the inequalities is strict) such that
$$
 \rho(\tilde H)\leq\frac{a+(a+1)(b-2)+\min\{k-1,a+2\}}{a+b}, \text{ while }
 \rho(H^{k-1}_{\ell})=\frac{k-2+(k-1)(\ell+1)}{k+\ell}.
$$
Let us prove that $\rho(\tilde H)<\rho(H^{k-1}_{\ell})$ for $\ell$ large enough.

First, let $a<k-2$. Then $\rho(\tilde H)\leq\frac{a+(a+1)(b-2)+a+2}{a+b}=\frac{a+1}{1+a/b}$. This bound increases as $b$ increases. Therefore, $\rho(\tilde H)\leq \frac{a+1}{1+a/(\ell+2)}\to a+1$ as $\ell\to\infty$. In contrast, $\rho(H^{k-1}_{\ell})\to k-1>a+2$ as $\ell\to\infty$. Therefore, $\rho(\tilde H)<\rho(H^{k-1}_{\ell})$ for $\ell$ large enough.

Finally, let $a=k-2$. Then 
$$
\rho(\tilde H)\leq\frac{k-2+(k-1)(b-1)}{k-2+b}=\frac{k-1}{1+(k-2)/b}-\frac{1}{k-2+b}=:\zeta(b).
$$
Since $\zeta(b)$ increases, the bound is strictly less than $\zeta(\ell+2)=\rho(H^{k-1}_{\ell})$. $\Box$\\

It is clear that the following sentence in $\mathcal{L}^k$ expresses the property `having $k-2$ vertices $v_1,\ldots,v_{k-2}$ such that the minimum length of a terminated $(k-1)$-chain rooted in $v_1,\ldots,v_{k-2}$ equals $\ell$':
$$
 \varphi_{\ell}^{k-1}=\exists v_1\ldots\exists v_{k-2}\exists x_0\,\left(\left[\bigwedge_{i\neq j}(\neg [v_i=v_j] \wedge \neg [v_i\sim v_j] ) \right]\wedge\left[\bigwedge_{i=1}^{k-2}x_0\sim v_i\right]\wedge
 {\sf TC}_{\ell}\wedge\left[\bigwedge_{i=1}^{\ell-1}\neg {\sf TC}_i\right]\right),
$$
where
$$
 {\sf TC}_i={\sf TC}_i(v_1,\ldots,v_{k-2},x_0)=
 \exists x_1\,\, \biggl(N_{x_1}(x_0,v_1,\ldots,v_{k-2})\wedge
$$
$$ 
 \biggl[\exists x_0\,\,\biggl(N_{x_0}(x_1,v_1,\ldots,v_{k-2})\wedge\ldots \biggl[\exists v_{k-2}\,\,N_{v_{k-2}}(x_0,x_1,v_1,\ldots,v_{k-3})\biggr]\ldots\biggr)\biggr]\biggr),
$$
$$
N_a(b_1,\ldots,b_{k-1})=\left(\bigwedge_{i=1}^{k-1} a\sim b_i\right).
$$
Clearly, if, in the sequence of vertices $x_1,x_0,x_1,\ldots$ whose existence is stated by ${\sf TC}_i$, there is a pair of coincident vertices, then this contradict with $\neg{\sf TC}_i$ for some $i\in\{1,\ldots,\ell-1\}$. The last vertex $v_{k-2}$ should not be equal to the vertex denoted earlier by the same variable, since the former $v_{k-2}$ should not be adjacent to any of $v_1,\ldots,v_{k-3}$. It is easy to see that, if this last $v_{k-2}$ coincides with one of the vertices in the sequence $x_1,x_0,x_1,\ldots$, this contradicts with one of $\neg{\sf TC}_i$ as well. Indeed, consider non-adjacent vertices $v_1,\ldots,v_{k-2}$, their common neighbor $x_0$, and the vertices $x_1,\ldots,x_{\ell}$ such that, for every $i\in\{1,\ldots,\ell\}$, $x_i$ is adjacent to $x_{i-1}$ and to $v_1,\ldots,v_{k-2}$. Assume that, for some $j\in\{0,1,\ldots,\ell-2\}$, the vertex $x_j$ is adjacent to both $x_{\ell-1}$ and $x_{\ell}$. Then the vertices $v_1,\ldots,v_{k-2},x_0,x_1,\ldots,x_j,x_{\ell-1},x_{\ell}$ form a terminated $1/(k-1)$-chain of length $j+1<\ell$.\\

It remains to prove that ${\sf P}(G(n,n^{-1/(k-1)})\models\phi)$ does not converge neither to $0$ nor to $1$. Since $\rho(H^{k-1}_{\ell})=\frac{k-2+(k-1)(\ell+1)}{k+\ell}$, it converges to $k-1$ and smaller than $k-1$ for all $\ell$. Therefore, for every $\varepsilon>0$, there exists $\ell^0$ such that, for all $\ell\geq\ell^0$, $\rho(H^{k-1}_{\ell})\in(\frac{1}{\varepsilon+1/(k-1)},k-1).$ Fix $\ell\geq\max\{\ell_0,\ell^0\}$ (where $\ell_0$ is from Lemma~\ref{H_str_b}). Then $H^{k-1}_{\ell}$ is strictly balanced, and 
$$
\frac{1}{\rho(H^{k-1}_{\ell})}\in\left(\frac{1}{k-1},\frac{1}{k-1}+\varepsilon\right).
$$
Set $\alpha=\frac{1}{\rho(H^{k-1}_{\ell})}$.
By Theorem~\ref{pois}, ${\sf P}(G(n,n^{-\alpha})\supset H^{k-1}_{\ell})\to 1-e^{-1/a(H^{k-1}_{\ell})}\in(0,1)$ as $n\to\infty$. It is clear that $\varphi_{\ell}^{k-1}$ implies existence of a subgraph isomorphic to $H^{k-1}_{\ell}$, and so 
\begin{equation}
\overline{\lim}_{n\to\infty}{\sf P}\left(G(n,n^{-\alpha})\models\varphi_{\ell}^{k-1}\right)\leq 1-e^{-1/a(H^{k-1}_{\ell})}.
\label{upper}
\end{equation}
Nevertheless, if a graph $G$ contains a subgraph isomorphic to $H^{k-1}_{\ell}$, but does not contain any subgraph that contains $H^{k-1}_{\ell}$ as a spanning proper subgraph, and does not contain a copy $X$ of $H^{k-1}_{\ell}$ and a vertex outside $X$ having at least $k-1$ neighbors in $V(X)$, then $G\models\varphi_{\ell}^{k-1}$. Clearly, if $H$ contains $H^{k-1}_{\ell}$ as a spanning proper subgraph, then $\rho(H)>1/\alpha$. By Theorem~\ref{small_thr}, a.a.s. in $G(n,n^{-\alpha})$, there is no copy of such $H$. Finally, if $H$ is obtained by adding to $H^{k-1}_{\ell}$ a vertex with at least $k-1$ neighbors, then $\rho(H)>1/\alpha$ as well, and so, a.a.s. in $G(n,n^{-\alpha})$, there is no copy of such $H$. This immediately implies that
\begin{equation}
\underline{\lim}_{n\to\infty}{\sf P}\left(G(n,n^{-\alpha})\models\varphi_{\ell}^{k-1}\right)\geq 1-e^{-1/a(H^{k-1}_{\ell})}.
\label{lower}
\end{equation}
It follows from~(\ref{upper})~and~(\ref{lower}) that $\lim_{n\to\infty}{\sf P}\left(G(n,n^{-\alpha})\models\varphi_{\ell}^{k-1}\right)=1-e^{-1/a(H^{k-1}_{\ell})}$, and this finishes the proof.

\section{Acknowledgements}
The paper is supported by the grant 18-71-00069 of Russian Science Foundation.

\end{document}